\documentclass[12pt]{article}
\usepackage{bm,hyperref}

\begin{document}
\title{Connections between Romanovski and other polynomials}
\author{H. J. Weber\\Department of Physics\\University of Virginia\\
Charlottesville, VA 22904, USA}

\maketitle
\begin{abstract}
A connection between Romanovski polynomials and those polynomials that solve 
the one-dimensional Schr\"odinger equation with the trigonometric Rosen-Morse 
and hyperbolic Scarf potential is established. The map is constructed by 
reworking the Rodrigues formula in an elementary and natural way. The 
generating function is summed in closed form from which recursion relations 
and addition theorems follow. Relations to some classical polynomials are also 
given.    
\end{abstract}
\vspace{3ex}
\leftline{MSC: 33C45, 42C15, 33C30, 34B24}
\leftline{PACS codes: 02.30.Gp; 02.30.Hq; 02.30.Jr, 03.65.Ge}
\leftline{Keywords: Romanovski polynomials; complexified Jacobi polynomials;}
\leftline{~~~~~~~~~~~~~~~generating function; recursion relations; addition 
theorems}

\section{Introduction and review of basic results}

Romanovski polynomials were discovered in 1884 by Routh~\cite{routh} in the 
form of complexified Jacobi polynomials on the unit circle in the complex 
plane and were then rediscovered as real polynomials by Romanovski~\cite{rom} 
in a statistics framework. Recently real polynomial solutions of the 
Scarf~\cite{cot} and Rosen-Morse potentials~\cite{ck} in (supersymmetric) 
quantum mechanics were recognized~\cite{rwack} to be related to the Romanovski 
polynomials.       

Here we apply to Romanovski polynomials a recently introduced natural method 
of reworking the Rodrigues formula~\cite{hjw} that leads to connections 
with other polynomials. 
 
The paper is organized as follows. In the Introduction we define the 
complementary polynomials $Q_\nu^{(\alpha,-a)}(x)$ then establish both a 
recursive differential equation satisfied by them and the procedure for the 
systematic construction of the $Q_{\nu}^{(\alpha,-a)}(x),$  derive their 
Sturm-Liouville differential equation (ODE), their generating function and its 
consequences, all based on the results of ref.~\cite{hjw}. Section~2 deals 
with a parameter addition theorem, Sect.~3 with orthogonality integrals and 
Sect.~4 with connections of these Romanovski polynomials to the classical 
polynomials. Sect.~5 deals with further applications using auxiliary 
polynomials.  

{\bf Definition.} The Rodrigues formula that generates the polynomials is 
given by 
\begin{eqnarray}
P_l^{(a,\alpha)}(x)=\frac{1}{w_l(x)}\frac{d^{l}}{dx^{l}}[w_l(x)\sigma(x)^{l}]
=\frac{\sigma(x)^l}{w_0(x)}\frac{d^l w_0(x)}{dx^l},
 \quad l=0,1,...\ .
\label{rds}
\end{eqnarray}
where $\sigma(x)\equiv 1+x^2$ is the coefficient of $y''$ of the hypergeometric
 ODE (1) of ref.~\cite{hjw} that the polynomials satisfy. The variable $x$ is 
real and ranges from $-\infty$ to $+\infty.$ The corresponding weight functions
\begin{eqnarray}
w_l(x)=\sigma(x)^{-(l+a+1)} e^{-\alpha \cot^{-1} x}=\sigma(x)^{-l} w_0(x)
\label{wf}
\end{eqnarray}
depend on the parameters $a, \alpha$ that are independent of the degree $l$ of 
the polynomial $P_l^{(a,\alpha)}(x).~\diamond$ 

{\it Lemma.} {\it The weight functions of the $P_l^{(a,\alpha)}(x)$ polynomials
 satisfy Pearson's first-order ODE} 
\begin{eqnarray}
\sigma(x)\frac{dw_l(x)}{dx}=[\alpha-2x(l+a+1)] w_l(x).  
\label{pearson}
\end{eqnarray} 
{\bf Proof.} This is straightforward to check using $d\cot^{-1} x/dx=-1/
\sigma(x).~\diamond$
 
We now apply the method of ref.~\cite{hjw} and introduce the {\it complementary
 polynomials} $Q_\nu^{(\alpha,-a)}(x)$ defining them inductively according to 
the Rodrigues representation (see Eq.~(5) of ref.~\cite{hjw}) 
\begin{equation}
P_l^{(a,\alpha)}(x)=\frac{1}{w_l(x)}\frac{d^{l-\nu}}{dx^{l-\nu}}
[\sigma(x)^{l-\nu}w_l(x)Q_\nu^{(\alpha,-a)}(x)],\,\,\,\,\, \nu=0, 1, \ldots, l.
~\diamond
\label{Qdef}
\end{equation}
For Eq.~(\ref{Qdef}) to agree with the Rodigues formula~(\ref{rds}) for $\nu=0$
 requires $Q_0^{(\alpha,-a)}(x)\equiv 1.$ Then comparing Eq.~(\ref{pearson}) 
with Eq.~(4) of ref.~\cite{hjw} we find the coefficient $\tau(x)=\alpha-
\sigma'(x)(a+l)$ of $y'$ in the ODE of the polynomials. Comparing instead with 
the ODE~(37) of ref.~\cite{rwack} gives their parameter $\beta=-a.$ 
 
We now {\it identify the polynomials of ref.~\cite{hjw} with those defined in  
Eq.~(\ref{Qdef})}  
\begin{eqnarray}
{\cal P}_\nu(x;l)=Q_\nu^{(\alpha,-a)}(x),~l\geq \nu.\label{id}
\end{eqnarray}
We will show below that the polynomials $Q_\nu^{(\alpha,-a)}(x)$ are 
independent of the parameter $l.~\diamond$

{\it Theorem~1.1. $Q_\nu^{(\alpha,-a)}(x)$ is a polynomial of degree $\nu$ 
that satisfies the recursive differential relation}  
\begin{eqnarray}\nonumber
Q_{\nu+1}^{(\alpha,-a)}(x)&=&\sigma(x)\frac{dQ_{\nu }^{(\alpha,-a)}(x)}{dx}
+[\tau(x)+2x(l-\nu-1)]Q_{\nu}^{(\alpha,-a)}(x)\\\nonumber
&=&\sigma(x)\frac{dQ_{\nu }^{(\alpha,-a)}(x)}{dx}
+[\alpha-2x(a+\nu+1)]Q_{\nu}^{(\alpha,-a)}(x),\\
\nu&=&0, 1, \ldots.\  
\label{rode}
\end{eqnarray}  
{\bf Proof.} The inductive proof of Theorem~2.2 of ref.~\cite{hjw} applied to 
the polynomial $Q_\nu^{(\alpha,-a)}(x)$ proves this theorem, and 
Eq.~(\ref{rode}) agrees with Eq.~(76) of ref.~\cite{rwack} provided their 
parameter $\beta=-a$ in ref.~\cite{hjw}. Since Eq.~(\ref{rode}) is independent 
of the parameter $l,$ so are the polynomials $Q_\nu^{(\alpha,-a)}(x)$ that are 
generated from it. $\diamond$ 

Comparing the recursive ODE~(\ref{rode}) with one stated in~\cite{rom} leads 
us to the {\it identification of our polynomials}  
\begin{eqnarray}
Q_{k}^{(\alpha,k-m)}(x)=\varphi_k(m,x),\ Q_\nu^{(\alpha,-a)}(x)=
\varphi_{\nu}(a+\nu, x)
\label{roq}
\end{eqnarray}
as a {\it Romanovski polynomial} (with its parameter depending on its degree), 
and comparing with Eq.~(69) of ref.~\cite{rwack},
\begin{eqnarray}
Q_\nu^{(\alpha,-a)}(x)=R_\nu^{(\alpha,\beta-\nu)}(x),~\beta=-a.~\diamond 
\label{comp}
\end{eqnarray}
Notice that the parameter $\alpha$ ($-\nu$ in ref.~\cite{rom}) is suppressed 
in Romanovski's notation. The fact that the integer index $\nu$ of the 
complementary polynomials occurs in the parameter $m$ of the Romanovski 
polynomials is sometimes disadvantageous (for orthogonality), but also 
occasionally a definite advantage (for the generating function). Moreover,
\begin{eqnarray}
{\cal P}_l(x;l)=Q_l^{(\alpha,-a)}(x)=K_l C_l^{(\alpha,-a)}(x),
\label{cs}
\end{eqnarray}
where $K_l$ is a normalization constant and the $C_l^{(\alpha,-a)}(x),$ after 
a change of variables, become part of the solutions of the Scarf and 
Rosen-Morse potentials in the Schr\"odinger equation~\cite{cot},\cite{ck}. 
However, for $C_l^{(\alpha,-a)}(x)$ to be part of the solution of the 
Schr\"odinger equation with the trigonometric Rosen-Morse potential requires 
$\alpha=\alpha_l=\frac{2b}{l+a}$~\cite{ck}; but the polynomials may be defined 
for the general parameter $\alpha.$ 

Recursion relations and recursive ODEs are practical tools to systematically 
generate the polynomials. 

{\it Theorem~1.2. The polynomial $Q_\nu^{(\alpha,-a)}(x)$ satisfies the basic 
recursive ODE}
\begin{eqnarray}
\frac{dQ^{(\alpha,-a)}_{\nu}(x)}{dx}=-\nu(2a+\nu+1)Q^{(\alpha,-a)}_{\nu-1}(x)
\equiv -\lambda_\nu Q^{(\alpha,-a)}_{\nu-1}(x).
\label{bode}
\end{eqnarray}
{\bf Proof.} Eq.~(\ref{bode}) follows from a comparison of the recursive 
ODE~(\ref{rode}) with a three-term recursion relation as outlined in 
Corollary~4.2 of ref.~\cite{hjw}, is ODE~(32) of ref.~\cite{hjw},  and agrees 
with Eq.~(75) of ref.~\cite{rwack} provided their $\beta=-a,$ which is 
consistent with our previous statements. $\diamond$

Thus, taking a derivative of $Q^{(\alpha,-a)}_{\nu}(x)$ just lowers its degree 
(and index) by unity, up to a constant factor, a property the Romanovski 
polynomials share with all classical polynomials.   

{\it Theorem~1.3. The polynomials $Q_{\nu}^{(\alpha,-a)}(x)$ satisfy the 
differential equation of Sturm-Liouville type}
\begin{eqnarray}
\sigma(x)\frac{d^2 Q^{(\alpha,-a)}_\nu(x)}{dx^2}+[\alpha-\sigma'(x)(a+\nu)]
\frac{d Q^{(\alpha,-a)}_\nu(x)}{dx}=-\lambda_\nu Q^{(\alpha,-a)}_\nu(x).
\label{qode}
\end{eqnarray}
{\bf Proof.} Substituting the basic ODE~(\ref{bode}) in the recursive 
ODE~(\ref{rode}) yields
\begin{eqnarray}
Q^{(\alpha,-a)}_{\nu+1}(x)=-\frac{\sigma(x)}{\lambda_{\nu+1}}
\frac{d^2Q^{(\alpha,-a)}_{\nu+1}(x)}{dx^2}-\frac{1}{\lambda_{\nu+1}}[\alpha-
(a+\nu+1)\sigma']\frac{dQ^{(\alpha,-a)}_{\nu+1}(x)}{dx}
\end{eqnarray}  
which, for $\nu\to \nu-1$ is the ODE of the theorem. Again, the ODE is 
independent of the parameter $l$ in Eq.~(\ref{id}). $\diamond$ 

{\it Theorem~1.4. The polynomial $P_l(x)$ satisfies the ODE} 
\begin{eqnarray}
\sigma(x)\frac{d^2 P_l(x)}{dx^2}+\tau(x)\frac{d P_l(x)}{dx}=-\lambda_l P_l(x).
\label{pode}
\end{eqnarray} 
{\bf Proof.}  For $\nu=l,$ we use $P_l(x)={\cal P}_l(x;l)$ in the notation of 
ref.~\cite{hjw} to rewrite the recursive ODE~(\ref{rode}) of ref.~\cite{hjw} 
and Eq.~(\ref{bode}) as 
\begin{eqnarray}\nonumber
{\cal P}_l(x;l)&=&\sigma(x){\cal P'}_{l-1}(x;l)+\tau(x){\cal P}_{l-1}(x;l)
=P_l(x)\\&=&-\frac{\sigma(x)}{\lambda_l}P''_l(x)-\frac{\tau(x)}{\lambda_l}
P'_l(x),
\label{Qode}
\end{eqnarray}
which is the ODE~(1) in ref.~\cite{hjw} for the polynomial $P_l(x).$ (Note 
that $\tau(x)$ is given right after Eq.~(\ref{Qdef}).) $\diamond$ 

{\it Theorem~1.5. The polynomials $Q^{(\alpha,-a)}_{\nu}(x)$ satisfy the 
generalized Rodrigues formulas}
\begin{eqnarray}
Q^{(\alpha,-a)}_{\nu}(x)&=&w^{-1}_l(x) \sigma(x)^{\nu-l}
\frac{d^{\nu}}{dx^{\nu}}[w_l(x)\sigma(x)^{l}]=\frac{\sigma(x)^\nu}{w_0(x)}
\frac{d^{\nu}w_0(x)}{dx^{\nu}};
\label{qrod}\\\nonumber 
Q^{(\alpha,-a)}_{\nu}(x)&=&w_l^{-1}(x)\sigma(x)^{\nu-l}\frac{d^{\nu-\mu}}
{dx^{\nu-\mu}}\left(\sigma(x)^{l-\mu}w_l(x) Q^{(\alpha,-a)}_{\mu}(x)\right),
\\ \mu&=&0, 1, \ldots, \nu .
\label{gqrod}
\end{eqnarray}
{\bf Proof.} These Rodrigues formulas are those of Theorem~2.3 of 
ref.~\cite{hjw}; they agree with Eqs.~(72) and (73) of ref.~\cite{rwack} 
provided their $\beta=-a,$ as we found earlier. Note that, from Eq.~(\ref{wf}) 
the product $w_l(x)\sigma(x)^l$ does not depend on $l,$ so there is no $l$ 
dependence in Eqs.~(\ref{qrod},\ref{gqrod}). $\diamond$

The $Q^{(\alpha,-a)}_n(x)$ polynomials generalize the $P_n(x)$ in the sense of 
allowing any power of $\sigma(x)$ in the Rodrigues formula, not just 
$\sigma(x)^n$ as for the $P_n(x).$ In other words, the $Q^{(\alpha,-a)}_n(x)$ 
are associated $P_n(x)$ (or $C_n^{(\alpha,-a)}(x))$ polynomials, as in the 
relationship between Laguerre (Legendre) and associated Laguerre (Legendre) 
polynomials. 

A generalization of Eq.~(\ref{qrod}), 
\begin{eqnarray}
Q^{(\alpha,-a-l)}_{\nu}(x)=\frac{\sigma(x)^{\nu+l}}{w_0(x)}
\frac{d^{\nu}}{dx^{\nu}}\left(\sigma(x)^{-l}w_0(x)\right),~l=0,\pm 1,\ldots,
\end{eqnarray}
just reproduces the same polynomial with a shifted parameter $a\to a+l.$ 

{\it The generating function for the $Q_\nu^{(\alpha,-a)}(x)$ polynomials is 
defined as}   
\begin{eqnarray}
Q(x,y;\alpha,-a)=\sum_{\nu=0}^{\infty}\frac{y^{\nu}}{\nu!} 
Q^{(\alpha,-a)}_{\nu}(x).~\diamond 
\label{genf}
\end{eqnarray}

The generating function is our main tool for deriving recursion relations.
 
{\it Theorem~1.6 The generating function can be summed in the closed form}   
\begin{eqnarray}\nonumber
w_l(x)Q(x,y;\alpha,-a)&=&\sigma(x)^{1-l}\sum_{\nu=0}^{\infty}
\frac{(y\sigma(x))^{\nu}}{\nu!}\frac{d^{\nu}}{dx^{\nu}}
\left(\sigma(x)^{-(a+1)}e^{\alpha\cot^{-1}x}\right)\,,\\\nonumber
Q(x,y;\alpha,-a)&=&(1+x^2)^{a+1}e^{\alpha\cot^{-1}x}
[1+(x+y(1+x^2))^2]^{-(a+1)}\\&\cdot&e^{-\alpha\cot^{-1}(x+y(1+x^2))}\ ,
\label{3_8}
\end{eqnarray}
\begin{eqnarray}\nonumber
&&\frac{\partial^{\mu}}{\partial y^{\mu}}\left( w_l(x)\sigma(x)^{l-1}
Q(x,y;\alpha,-a)\right)\\&=&\sigma(x)^{\mu}\sum_{\nu=\mu}^{\infty}
\frac{(y\sigma(x))^{\nu-\mu}}{(\nu-\mu)!}\frac{d^{\nu-\mu}}{dx^{\nu-\mu}}
(\sigma(x)^{-(a+\mu)}e^{-\alpha\cot^{-1}x}Q^{(\alpha,-a)}_{\mu}(x)),\\\nonumber
&&\frac{\partial^{\mu}Q(x,y;\alpha,-a)}{\partial y^{\mu}}=w_l(x)^{-1}
\sigma(x)^{\mu-l}[(1+(x+y\sigma(x))^2]^{-(a+\mu)}\\ &\cdot& 
e^{-\alpha\cot^{-1}(x+y\sigma(x))^2}Q^{(\alpha,-a)}_{\mu}((x+y\sigma(x))^2) . 
\label{3_9}
\end{eqnarray}
{\it Both Taylor series converge if $x$ and $x+y\sigma(x)$ are regular points 
of the weight function.}

{\bf Proof.} The first relation is derived in Theorem~3.2 of ref.~\cite{hjw} 
by substituting the Rodrigues formula~(\ref{qrod}) in the defining 
series~(\ref{genf}) of the generating function and recognizing it as a Taylor 
series. The other follows similarly. $\diamond$

{\it Theorem~1.7. The generating function satisfies the partial differential 
equation (PDE)}
\begin{eqnarray}
\frac{\partial Q(x,y;\alpha,-a)}{\partial y}=
\frac{\sigma(x)Q(x,y;\alpha,-a)}{1+[x+y\sigma(x)]^2}
Q^{(\alpha,-a)}_1(x+y\sigma(x)). 
\label{pde}
\end{eqnarray} 
{\bf Proof.} This PDE is derived by straightforward differentiation in 
Theorem~3.3 of ref.~\cite{hjw} in preparation for recursion relations by 
translating the case $\mu=1$ in Eq.~(\ref{3_9}) into a partial differential 
equation (PDE). $\diamond$ 

One of the main consequences of Theorem~1.7 is a general recursion relation.
  
{\it Theorem~1.8. The $Q^{(\alpha,-a)}_{\nu}(x)$ polynomials satisfy the
three-term recursion relation}
\begin{eqnarray}\nonumber
Q^{(\alpha,-a)}_{\nu+1}(x)=[\alpha-2x(a+\nu+1)]Q^{(\alpha,-a)}_{\nu}(x)
-\nu\sigma(x)(2a+\nu+1)Q^{(\alpha,-a)}_{\nu-1}(x).\\
\label{recu3}
\end{eqnarray}

{\bf Proof.} Equation~(\ref{pde}) translates into this recursion relation by 
substituting Eq.~(\ref{genf}) defining the generating function and thus 
rewriting this as
\begin{eqnarray}\nonumber
&&(1+y^2\sigma^2(x)+2xy)\sum_{\nu=1}^{\infty}\frac{y^{\nu-1}}{(\nu-1)!}
Q^{(\alpha,-a)}_{\nu}(x)\\&=&[\alpha-2x(a+1)-2y(a+1)\sigma(x)]
\sum_{\nu=0}^{\infty}\frac{y^{\nu}}{\nu!}Q^{(\alpha,-a)}_{\nu}(x).~\diamond
\end{eqnarray}

Just like the recursive ODE~(\ref{rode}), this recursion allows for a 
systematic construction of the Romanovski polynomials, in contrast to the 
Rodrigues formulas which become impractical for large values of the degree 
$\nu.$

{\it Theorem~1.9. The polynomials $Q_{\nu}^{(\alpha,-a)}(x)$ satisfy the 
differential equation of Sturm-Liouville type}  
 \begin{eqnarray}\nonumber
\frac{d}{dx}\left(\sigma(x)^{l-\nu+1}w_l(x)\frac{dQ_{\nu}^{(\alpha,-a)}(x)}
{dx}\right)&=&-\lambda_{\nu}\sigma(x)^{l-\nu}(x)w_l(x)  
Q_{\nu}^{(\alpha,-a)}(x);\\\lambda_{\nu}&=&\nu(2a+\nu+1)\ . 
\label{qsl}
\end{eqnarray}
{\bf Proof.} This ODE is equivalent to the ODE~(\ref{qode}) and agrees with 
Eq.~(78) of ref.~\cite{rwack} if $\beta=-a$ there. Note that the inductive 
proof in Theorem~5.1 in ref.~\cite{hjw} is much lengthier than our proof of 
Eqs.~(\ref{qode},\ref{pode}). $\diamond$ 

\section{Parameter Addition}

The multiplicative structure of the generating function of Eq.~(\ref{3_8})
involving the two parameters in the exponents of two separate functions, 
as displayed in 
\begin{eqnarray}
Q\left(x,\frac{y-x}{\sigma(x)};\alpha,-a\right)=\left(\frac{\sigma(x)}
{\sigma(y)}\right)^{a+1}e^{\alpha(\cot^{-1}x-\cot^{-1}y)},
\end{eqnarray}
allows for the following theorems. 

{\it  Theorem~2.1. The $Q_{\nu}^{(\alpha,-a)}(x)$ polynomials 
satisfy the parameter addition relation}  
\begin{eqnarray}
Q_{N}^{(\alpha_1+\alpha_2,-a_1-a_2-1)}(x)=\sum_{\nu_1=0}^N
\left(N\atop \nu_1\right)Q_{\nu_1}^{(\alpha_1,-a_1)}(x)
Q_{N-\nu_1}^{(\alpha_2,-a_2)}(x).
\end{eqnarray}
{\bf Proof.} This formula follows from the Taylor expansion 
\begin{eqnarray}
\sum_{\nu_1, \nu_2=0}^{\infty}\frac{y^{\nu_1+\nu_2}}{\nu_1!\nu_2!}
Q_{\nu_1}^{(\alpha_1,-a_1)}(x)Q_{\nu_2}^{(\alpha_2,-a_2)}(x)=
\sum_{N=0}^{\infty}\frac{y^{N}}{N!}
Q_{N}^{(\alpha_1+\alpha_2,-a_1-a_2-1)}(x). 
\end{eqnarray}
of the {{\it generating function identity}
\begin{eqnarray}
Q(x,y;\alpha_1,-a_1)Q(x,y;\alpha_2,-a_2)=Q(x,y;\alpha_1+\alpha_2,-a_1-a_2-1).~
\diamond
\end{eqnarray}

Given the complexity of the polynomials, the elegance and simplicity of this 
relation are remarkable.
 
{\bf Example.} The case $N=0$ is trivial, and $N=1$ becomes the additive 
identity
\begin{eqnarray}\nonumber
Q_1^{(\alpha_1+\alpha_2,-(a_1+a_2+1))}(x)&=&\alpha_1+\alpha_2-2x(a_1+a_2+2)\\
\nonumber&=&[\alpha_1-2x(a_1+1)]+[\alpha_2-2x(a_2+1)]\\&=&
Q_1^{(\alpha_1,-a_1)}(x)+Q_1^{(\alpha_2,-a_2)}(x).
\end{eqnarray}
The first case involving additive and multiplicative aspects of the polynomials
 is $N=2$ which we decompose and multiply out as follows:
\begin{eqnarray}\nonumber
Q_{2}^{(\alpha_1+\alpha_2,-a_1-a_2-1)}&=&
[\alpha_1+\alpha_2-2x(a_1+a_2+2)][\alpha_1+\alpha_2-2x(a_1+a_2+3)]\\\nonumber 
&-&2\sigma(x)(a_1+a_2+2)=\{[\alpha_1-2x(a_1+1)]\\\nonumber &+&
[\alpha_2-2x(a_2+1)]\}\{[\alpha_1-2x(a_1+2)]\\\nonumber&+&[\alpha_2-2x(a_2+1)]
\}-2\sigma(x)[(a_1+1)+(a_2+1)]\\\nonumber&=&[\alpha_1-2x(a_1+1)][\alpha_1
-2x(a_1+2)]-2\sigma(x)(a_1+1)\\\nonumber
&+&[\alpha_2-2x(a_2+1)][\alpha_2-2x(a_2+2)]-2\sigma(x)(a_2+1)\\\nonumber
&+&[\alpha_1-2x(a_1+1)][\alpha_2-2x(a_2+1)]\\\nonumber&+&[\alpha_2-2x(a_2+1)]
[\alpha_1-2x(a_1+1)]\\\nonumber&=&
Q_{2}^{(\alpha_1,-a_1)}+Q_{2}^{(\alpha_2,-a_2)}+Q_{1}^{(\alpha_1,-a_1)}
Q_{1}^{(\alpha_2,-a_2)}\\&+&Q_{1}^{(\alpha_2,-a_2)}Q_{1}^{(\alpha_1,-a_1)}.~\diamond
\end{eqnarray}

The addition theorem is consistent with {\it the homogeneous polynomial 
theorem in the variables} $x, \alpha, \sqrt{\sigma}$ (without using 
$\sigma=x^2+1$) {\it which the polynomials satisfy} and can be generalized to 
an arbitrary number of parameters.
 
{\it Theorem~2.2. The $Q\nu^{(\alpha,-a)}(x)$ polynomials satisfy the more 
general polynomial identity} 
\begin{eqnarray}\nonumber
&&Q_N^{(\alpha_1+\alpha_2+\cdots +\alpha_n,-(a_1+a_2+\cdots +a_n+n-1))}(x)\\&=&
\sum_{0\leq \nu_j\leq N, \nu_1+\cdots +\nu_n=N+n}\frac{N!}
{\prod_1^n \nu_j!}\prod_{j=1}^n Q_{\nu_j}^{(\alpha_j,-a_j)}(x).
\end{eqnarray}

{\bf Proof.} This follows similarly from the Taylor expansion of the {\it 
product identity of $n$ generating functions}
\begin{eqnarray}\nonumber
\prod_{j=1}^n Q(x,y;\alpha_j,-a_j)=Q(x,y;\alpha_1+\alpha_2+\cdots +\alpha_n,
-(a_1+a_2+\cdots +a_n+n-1)).\\ 
\end{eqnarray} 
As an application of the parameter addition theorem we now separate the two
parameters $a$ and $\alpha$ into two sets of simpler polynomials 
$Q_{\nu}^{(0,-a)}$ and $Q_{\mu}^{(\alpha,1)}.$  To this end, we expand the 
generating functions in the {\it identity}
\begin{eqnarray}
Q(x,y;\alpha,-a)=Q(x,y;0,-a)Q(x,y;\alpha,1)
\label{decomp}
\end{eqnarray}
in terms of their defining polynomials. This yields  

{\it Theorem~2.3. The $Q_\nu^{(\alpha,-a)}(x)$ polynomials satisfy the 
decomposition identity}
\begin{eqnarray}
Q^{(\alpha,-a)}_N(x)=\sum_{\nu=0}^N\left(N\atop \nu\right)
Q^{(0,-a)}_{\nu}(x)Q^{(\alpha,1)}_{N-\nu}(x). 
\label{simpl}
\end{eqnarray}
{\bf Proof.} This identity follows from a Taylor expansion of the generating 
function identity~(\ref{decomp}) in terms of sums of products of polynomials 
involving only one parameter each. $\diamond$
 
{\bf Definition.} The generating function
\begin{eqnarray}
e^{\alpha [\cot^{-1}x-\cot^{-1}(x+y\sigma(x))]}=\sum_{\nu=0}^{\infty}
\frac{y^{\nu}}{\nu!}Q^{(\alpha,1)}_{\nu}(x)
\end{eqnarray}
defines the second set of the polynomials, while the first one will be treated
 in detail below upon expanding the polynomials $Q^{(0,-a)}_\nu(x)$ as finite 
sums of Gegenbauer polynomials in Sect.~IV or finite power series in 
Eq.~(\ref{expl}). $\diamond$

We also note that $Q_\nu^{(0,-a)}(x)=K_\nu C_\nu^{(0,-a)}(x),$ so the latter 
also have a Gegenbauer polynomial expansion.  

{\it Theorem~2.4. The $Q_\nu^{(\alpha,-a)}(x)$ polynomials satisfy the
(parity) symmetry relation} 
\begin{eqnarray}
Q_{\nu}^{(-\alpha,-a)}(x)=(-1)^{\nu}Q_{\nu}^{(\alpha,-a)}(-x). 
\end{eqnarray}

{\bf Proof.} This relation derives from the {\it generating function identity}
\begin{eqnarray}
Q(-x,-y;-\alpha,-a)=Q(x,y;\alpha,-a)
\end{eqnarray}
which holds because $\alpha\cot^{-1}(x+y\sigma(x))$ in the generating function 
 Eq.~(\ref{genf}) stays invariant under $\alpha\to -\alpha, x\to -x, y\to 
-y,$ and $\sigma(-x)=\sigma(x).~\diamond$ 

\section{Orthogonality Integrals}

This section deals with an application of the generating function to integrals
that are relevant for studying the orthogonality of the polynomials.
 
{\bf Definition.} We define orthogonality integrals for the 
$Q_\nu^{(\alpha,-a)}(x)$ polynomials by~\cite{rwack}
\begin{eqnarray}\nonumber
O^{(a,\alpha)}_{\mu,\nu}&=&\int_{-\infty}^{\infty}dx
\frac{Q^{(\alpha,-a)}_{\mu}(x)Q^{(\alpha,-a)}_{\nu}(x)e^{-\alpha\cot^{-1}x}}
{\sigma(x)^{(\mu+\nu)/2+a+2}}=0,\\a&>&-3/2,~\mu+\nu~\rm{even},
\label{orthdef}
\end{eqnarray}
while for $\mu+\nu$ odd there needs to be an extra $\sqrt{\sigma(x)}$ in the 
numerator for the orthogonality integrals to vanish. $\diamond$ 

Thus, the $Q_\nu^{(\alpha,-a)}(x)$ polynomials form two infinite subsets, each 
with general orthogonality, but polynomials from different subsets are not 
mutually orthogonal. While displaying infinite orthogonality, this property 
falls short of the general orthogonality of all classical polynomials. The 
$Q_\nu^{(\alpha,-a)}(x)$ polynomials form a partition of the set of all 
Romanovski polynomials, as shown in Eq.~(\ref{comp}), with upper index 
dependent on the running index $\nu,$ though. The Romanovski polynomials 
$R_\nu^{(\alpha,\beta)}(x)$ with upper indices independent of the degree $\nu,$
 the running index, form another partition that has the finite orthogonality, 
as discussed in more detail in ref.~\cite{rwack}. The orthogonality of the 
$C_n^{(\alpha_n,-a)}(x)$ polynomials from the Schr\"odinger equation with 
$\alpha=\alpha_n$ as discussed below Eq.~(\ref{cs}), is yet another form of 
orthogonality similar to that of hydrogenic wave functions, which also differs 
from the mathematical orthogonality of associated Laguerre polynomials, the 
subject of Exercise~13.2.11 in ref.~\cite{aw}. The orthogonality integrals 
of Eq.~(\ref{orthdef}) suggest analyzing the following integral of the 
generating functions
\begin{eqnarray}\nonumber
I(y)&=&\int_{-\infty}^{\infty}\frac{dx}{\sigma(x)^{a+2}}\left(
Q(x,\frac{y}{\sqrt{\sigma}};\alpha,-a)\right)^2 e^{-\alpha \cot^{-1}x}\\
\nonumber &=&\sum_{\nu_1,\nu_2=0}^{\infty}\frac{y^{\nu_1+\nu_2}}{\nu_1!
\nu_2!}\int_{-\infty}^{\infty}dx\frac{Q^{(\alpha,-a)}_{\nu_1}(x)
Q^{(\alpha,-a)}_{\nu_2}(x)e^{-\alpha\cot^{-1}x}}
{\sigma(x)^{(\nu_1+\nu_2)/2+a+2}}\\&=&\sum_{\nu_1,\nu_2=0}^{\infty}
\frac{y^{\nu_1+\nu_2}}{\nu_1!\nu_2!}O^{(a,\alpha)}_{\nu_1, \nu_2}
\label{orthint}
\end{eqnarray}
which is written directly in terms of orthogonality integrals
$O^{(\alpha,-a)}_{\nu_1, \nu_2}$ defined in Eq.~(\ref{orthdef}). On the other 
hand, we can express the integral as
\begin{eqnarray}\nonumber
I(y)&=&\int_{-\infty}^{\infty}\frac{dx e^{-\alpha\cot^{-1}x+2\alpha\cot^{-1}x
-2\alpha\cot^{-1}(x+y\sqrt{\sigma})}}{\sigma(x)^{a+2}[1+y^2+\frac{2xy}
{\sqrt{\sigma}}]^{2(a+1)}}\\\nonumber&=&\int_0^{\infty}\frac{dx}{\sigma^{a+2}}
\frac{e^{\alpha\cot^{-1}x-2\alpha\cot^{-1}(x+y\sqrt{\sigma})}}
{[1+y^2+\frac{2xy}{\sqrt{\sigma}}]^{2(a+1)}}\\&+&\int_0^{\infty}\frac{dx}
{\sigma^{a+2}}\frac{e^{-\alpha\cot^{-1}x+2\alpha\cot^{-1}(x-y\sqrt{\sigma})}}
{[1+y^2-\frac{2xy}{\sqrt{\sigma}}]^{2(a+1)}},
\end{eqnarray}
which is manifestly not even in the variable $y.$ If the $Q^{(\alpha,-a)}_{\nu}$ 
polynomials were orthogonal, then the double sum in $I$ of Eq.~(\ref{orthint}) 
would collapse to a single sum over normalization integrals multiplied by even 
powers of $y,$ i.e. $I$ would be an even function of $y.$ This result shows 
that the $Q^{(\alpha,-a)}_{\nu}$ polynomials are not orthogonal in the 
conventional sense. In fact, the extra $\sqrt{\sigma}$ in the orthogonality 
integrals $O_{2\nu,2\mu+1}^{(a,\alpha)}$ is not built into the generating 
function. In other words, the fact that $I(y)\neq I(-y)$ indirectly confirms 
that the $Q^{(\alpha,-a)}_{\nu}$ polynomials have more complicated 
orthogonality properties than the Romanovski polynomials with parameters that 
are independent of the degree of the polynomial, as discussed in more detail 
in ref.~\cite{rwack}.

Let us next consider the special parameter $\alpha=0$ and analyze similarly the
 integral
\begin{eqnarray}\nonumber
I_0(y)&=&\int_{-\infty}^{\infty}\frac{dx}{\sigma(x)^{a+2}}
Q(x,\frac{y}{\sqrt{\sigma}};0,-a)=\sum_{\nu=0}^{\infty}\frac{y^{\nu}}{\nu!}
\int_{-\infty}^{\infty}\frac{dx Q^{(0,-a)}_{\nu}(x)}{\sigma^{\nu/2+2+a}}\\
\nonumber&=&\sum_{\nu=0}^{\infty}\frac{y^{\nu}}{\nu!}O^{(a, 0)}_{\nu, 0}
=\int_0^{\infty}\frac{dx}{\sigma^{a+2}}\frac{1}{[1+y^2+\frac{2xy}
{\sqrt{\sigma}}]^{a+1}}+\int_0^{\infty}\frac{dx}{\sigma^{a+2}}\frac{1}
{[1+y^2-\frac{2xy}{\sqrt{\sigma}}]^{a+1}},\\
\label{int0}
\end{eqnarray}
which is an even function of $y.$ If the $Q^{(0,-a)}_{\nu}(x)$ are orthogonal
to $Q^{(0,-a)}_{0}(x)=1$ then the sum in Eq.~(\ref{int0}) will collapse to its
first term and $I_0$ is a constant. It is quite a surprise that this actually
happens in the interval $-1\leq y\leq 1$ for all parameter values $a$ for
which the integral $I_0$ converges. For example, $I_0(y)=r(a)\pi=$const. with
a rational number $r(a)$ that depends on the exponent $a,$ where $r(0)=1/2,
r(1)=3/4, r(2)=5/8, r(3)=35/64, r(4)=3^2\cdot7/2^7,$ if $a$ is a non-negative
integer; in general $I_0(y)=\sqrt{\pi}\Gamma(a+3/2)/\Gamma(a+2).$ For $y>|1|, 
I_0(y)$ varies and deviates from the constant value. From the structure of the 
integral, this anomalous behavior of $I_0$ is rather unexpected. Since for 
$\alpha=0$ parity is conserved in the ODE~(\ref{qode}), the orthogonality 
integrals $O_{2\nu,0}^{(a,0)}$ are zero and $O_{2\nu,2\mu}^{(a,0)}=0$ more 
generally. Since it is shown in ref.~\cite{rwack} that $O_{2\nu,2\mu}^{(a,0)}$ 
for $\mu\neq \nu$ vanish, the $Q^{(0,-a)}_{\nu}(x)$ polynomials are orthogonal 
in the conventional sense. Since each $C_l^{(0,-a)}(x)$ is proportional to 
$Q^{(0,-a)}_{l}(x),$ the $C_l^{(0,-a)}(x)$ polynomials are orthogonal. This is 
confirmed by $I_0$ and its constancy in the interval $-1\leq y\leq 1$ is thus 
proved.

The restriction to parameter value $\alpha=0$ can be removed:
\begin{eqnarray}\nonumber
I_1(y)&=&\int_{-\infty}^{\infty}\frac{dx}{\sigma(x)^{a+2}}
Q(x,\frac{y}{\sqrt{\sigma}};\alpha,-a)=\sum_{\nu=0}^{\infty}\frac{y^{\nu}}
{\nu!}O^{(a, \alpha)}_{\nu, 0}\\\nonumber
&=&\int_{-\infty}^{\infty}\frac{dx}{\sigma^{a+2}}\frac{e^{-\alpha\cot^{-1}
(x+y\sqrt{\sigma})}}{[1+y^2+\frac{2xy}{\sqrt{\sigma}}]^{a+1}}\\
&=&\int_0^{\infty}\frac{dx}{\sigma^{a+2}}\frac{e^{-\alpha\cot^{-1}(x+y
\sqrt{\sigma})}}{[1+y^2+\frac{2xy}{\sqrt{\sigma}}]^{a+1}}+\int_0^{\infty}
\frac{dx}{\sigma^{a+2}}\frac{e^{-\alpha\cot^{-1}(x+y\sqrt{\sigma})}}
{[1+y^2-\frac{2xy}{\sqrt{\sigma}}]^{a+1}},
\label{int1}
\end{eqnarray}
which is neither even in $y$ nor independent of $y.$ Therefore, if we wish
to find the normalizations of the $Q_{\nu}^{(\alpha,-a)}$ polynomials we have
to split up the generating function into its even and odd parts in $y$ and
integrate them separately, each with the proper power of $\sigma(x)$ in the 
orthogonality integral.

\section{Relations with Gegenbauer Polynomials}

The relation of the Romanovski polynomials as complexified Jacobi polynomials 
on the unit circle in the complex plane is described in detail in 
ref.~\cite{rwack}. Therefore, we focus here on relations with Gegenbauer 
polynomials.

We start with the simplest case of parameter values $a=0=\alpha,$ which also
happens to be relevant for physics~\cite{ck}, to derive from the generating
function an expression for $Q_{\nu}^{(\alpha,-a)}(x)$ in terms of a finite sum
of Gegenbauer polynomials. For $a=0=\alpha,$ Eq.~(\ref{3_8}) takes the explicit
 form
\begin{eqnarray}
Q(x,y;0,0)=(1+x^2)[1+(x+y(1+x^2))^2]^{-1}=\frac{1}{1+2xy+y^2(1+x^2)}.
\label{g00}
\end{eqnarray}

{\it Theorem~4.1. The $Q_{m}^{(0, 0)}(x)$ polynomials have the expansion into 
Gegenbauer polynomials}
\begin{eqnarray}
Q_{m}^{(0, 0)}(x)=m!\sum_{n=0}^{[m/2]}(-1)^nx^{2n}C_{m-2n}^{(n+1)}(-x),\
m=0, 1,\ldots .
\end{eqnarray}
{\bf Proof.} For $|xy|/|1+2xy+y^2|<1,$ the generating function 
identity~(\ref{g00}) may be expanded as an absolutely converging geometric 
series 
\begin{eqnarray}
Q(x,y;0,0)=\sum_{n=0}^{\infty}\frac{(-x^2y^2)^n}{(1+2xy+y^2)^{n+1}}.
\end{eqnarray}
Substituting the generating function of Gegenbauer polynomials~\cite{aw}
$C_l^{(n+1)}(x)$,
\begin{eqnarray}
(1-2xy+y^2)^{-(n+1)}=\sum_{l=0}^{\infty}C_l^{(n+1)}(x)y^l,
\end{eqnarray}
we obtain the expansion 
\begin{eqnarray}\nonumber
Q(x,y;0,0)&=&\sum_{n=0}^{\infty}(-x^2y^2)^n
\sum_{l=0}^{\infty}C_l^{(n+1)}(-x)y^l\\&=&\sum_{m=0}^{\infty}y^m
\sum_{n=0}^{[m/2]}(-1)^nx^{2n}C_{m-2n}^{(n+1)}(-x),
\end{eqnarray}
where $m=l+2n$ was used upon interchanging the summations, with $[m/2]$ the 
integer part of $m/2$. On comparing with Eq.~(\ref{3_8}) defining the 
generating function $Q(y, x;\alpha,-a)$ we obtain the expansion of the 
$Q_{m}^{(0, 0)}(x)$ polynomials as a finite sum of Gegenbauer polynomials of 
Theorem~4.1. $\diamond$

Since $Q_{m}^{(0, 0)}(x)=K_{m}C_{m}^{(0, 0)}(x),$ this result is 
also valid for the $C_{m}^{(0, 0)}(x)$ polynomials. It can be generalized to 
parameter values $a\neq 0:$ 

{\it Theorem~4.2. The $Q_{N}^{(0,-a)}(x)$ polynomials have the Gegenbauer 
polynomial expansion}
\begin{eqnarray}
Q_{N}^{(0,-a)}(x)=N!\sum_{n=0}^{[N/2]}\left(-a-1\atop n\right)x^{2n}
C^{(n+a+1)}_{N-2n}(-x).
\end{eqnarray}

{\bf Proof.} This relation follows from expanding the generating function
\begin{eqnarray}\nonumber
Q(x,y;0,-a)&=&\left( \frac{\sigma(x)}{1+x^2+y^2\sigma^2(x)+2xy\sigma(x)}
\right)^{a+1}=\frac{1}{[1+2xy+y^2+x^2y^2]^{a+1}}\\\nonumber
&=&\sum_{n=0}^{\infty}\left(-a-1\atop n\right)\frac{(x^2y^2)^n}
{[1+2xy+y^2]^{n+a+1}}\\\nonumber&=&\sum_{n, l=0}^{\infty}\left(-a-1\atop n
\right)(x^2y^2)^nC^{(n+a+1)}_l(-x) y^l\\&=&\sum_{N=0}^{\infty} y^{N}
\sum_{n=0}^{[N/2]}\left(-a-1\atop n\right)x^{2n}C^{(n+a+1)}_{N-2n}(-x)
\end{eqnarray}
in terms of the binomial series and then again using the generating functions
of the Gegenbauer polynomials. $\diamond$

{\it Theorem~4.3. The $Q_{N}^{(\alpha,-a)}(x)$ polynomials have the general 
Gegenbauer polynomial expansion} 
\begin{eqnarray}
\frac{1}{N!}Q_{N}^{(\alpha,-a)}(x)=\sum_{\nu=0}^N\left(N\atop \nu\right)
Q^{(\alpha,1)}_{N-\nu}(x)\sum_{n=0}^{[N/2]}\left(-a-1\atop n\right)x^{2n}
C^{(n+a+1)}_{N-2n}(-x). 
\end{eqnarray}

{\bf Proof.} Substituting the expansion of Theorem~4.2 into Eq.(\ref{simpl}) 
in which the Gegenbauer polynomials depend only on the parameter $a$ while the 
$Q^{(\alpha,1)}_{\nu}(x)$ depend only on $\alpha$ yields the desired expansion.
 $\diamond$

The Gegenbauer polynomials are well-known generalizations of Legendre 
polynomials. The {\it hyperbolic Gegenbauer} ODE
\begin{eqnarray}
\sigma(x)y''-(2\lambda+1)xy'+\Lambda_l^{(\lambda)}y=0
\label{ggode}
\end{eqnarray}
becomes the ODE~(\ref{qode}) for $\alpha=0, \nu=l$ and $2\lambda+1=2(a+\nu),$ 
so the solutions of Eq.~(\ref{ggode}) are the $Q_{l}^{(0,-(\lambda-l+1/2))}(x)$
 polynomials. In fact, for $\alpha=0$ we can directly solve the 
ODE~(\ref{qode}) for the $Q_{\nu}^{(0,-a)}(x)$ polynomial solutions in terms 
of finite power series. 

{\it Theorem~4.4. The $Q_{N}^{(0,-a)}(x)$ polynomials have the explicit finite 
power series} 
\begin{eqnarray}\nonumber
Q_{N}^{(0,-a)}(x)&=&\sum_{\mu=0}^{[N/2]}x^{N-2\mu}a_{\mu},\
a_{\mu}=-\frac{(N-2\mu+2)(N-2\mu+1)}{2\mu(2a+2\mu+1)}a_{\mu-1},\\\nonumber
a_1&=&-\frac{N(N-1)}{2(2a+3)},\ a_{\mu}=(-1)^{\mu}
\frac{N(N-1)\cdots (N-2\mu+1)}{2^{\mu}\mu!(2a+3)(2a+1)\cdots
(2a+2\mu+1)}\\
\label{expl}
\end{eqnarray}
\begin{eqnarray}
Q_{N}^{(\alpha,-a)}(x)=\sum_{\nu=0}^N\left(N\atop \nu\right)
Q^{(\alpha,1)}_{N-\nu}(x)\sum_{\mu=0}^{[N/2]}x^{N-2\mu}a_{\mu}
\end{eqnarray} 
with $a_\mu$ from Eq.~(\ref{expl}) and $[N/2]$ denoting the integer part of 
$N/2.$

{\bf Proof.} Since the proof by mathematical induction is straightforward, we 
just give the results. As the ODE is invariant under the parity transformation,
 $x\to -x,$ we have even and odd solutions. Substituting Eq.~(\ref{expl}) in 
Eq.~(\ref{simpl}) yields the second relation stated in Theorem~4.4. $\diamond$

This is also valid for the $C_l^{(0,-a)}(x)$ polynomials up to the 
normalization $K_l.$ 

\section{Auxiliary polynomials}

Carrying out the innermost derivative of the Rodrigues formula~(\ref{rds}), we 
find  
\begin{eqnarray}\nonumber
P_l(x)&=&\frac{\sigma^l}{w_0}\frac{d^{l-1}}{dx^{l-1}}\left(\frac{w_0}{\sigma}
[\alpha-\sigma'(a+1)]\right)\\&=&\alpha Q^{(\alpha,-a-1)}_{l-1}(x)-
(a+1)\frac{\sigma^l}{w_0}\frac{d^{l-1}}{dx^{l-1}}\left(\frac{\sigma' w_0}
{\sigma}\right),
\label{aux}
\end{eqnarray}
and are led to define {\it auxiliary polynomials:}
\begin{eqnarray}
S_{l+1}(x)=\frac{\sigma(x)^{l+1}}{w_0(x)}\frac{d^l}{dx^l}
\left(\frac{\sigma'(x)w_0(x)}{\sigma(x)}\right).~\diamond
\label{defx}
\end{eqnarray}

{\bf Example.}
\begin{eqnarray}
S_1(x)=\sigma'(x),~S_2(x)=\sigma''\sigma(x)+\sigma'(x)[\alpha-\sigma'(x)(a+2)],
\ldots  .~\diamond
\end{eqnarray}
So
\begin{eqnarray}
S_l(x)=\frac{\alpha}{a+1}Q^{(\alpha,-a-1)}_{l-1}(x)-\frac{P_l(x)}{a+1}.
\label{alt1}
\end{eqnarray}
Applying a derivative to $w_0S_l/\sigma^l$ yields
\begin{eqnarray}
\frac{d^l}{dx^l}\left(\frac{\sigma'(x)w_0(x)}{\sigma(x)}\right)=\frac{\alpha}
{a+1}\frac{d}{dx}\left(\frac{w_0}{\sigma^l}Q^{(\alpha,-a-1)}_{l-1}(x)\right)
-\frac{1}{a+1}\frac{d}{dx}\left(\frac{w_0 P_l(x)}{\sigma^l}\right).
\end{eqnarray} 
Using the recursive ODEs for $Q^{(\alpha,-a-1)}_{l-1}$ and $P_l$ yields 
\begin{eqnarray}\nonumber
\frac{\sigma(x)^{l+1}}{w_0(x)}\frac{d^l}{dx^l}\left(\frac{\sigma'(x)w_0(x)}
{\sigma(x)}\right)&=&\frac{\alpha}{a+1}(Q^{(\alpha,-a-1)}_{l}(x)-\sigma'
Q^{(\alpha,-a-1)}_{l-1}(x))\\&-&\frac{Q^{(\alpha,-a)}_{l+1}(x)}{a+1}.
\label{alt2}
\end{eqnarray}
A comparison of Eqs.~(\ref{alt1},\ref{alt2}) yields 
\begin{eqnarray}
P_{l+1}(x)=\alpha\sigma'(x)Q^{(\alpha,-a-1)}_{l-1}(x)+Q^{(\alpha,-a)}_{l+1}(x),
\end{eqnarray}
complementing the relation ${\cal P}_l(x;l)=K_lC_l^{(\alpha,-a)}(x)
=Q^{(\alpha,-a)}_l(x)=P_l(x).$ 

For Laguerre polynomials $\sigma(x)=x$ and relation~(\ref{aux}) corresponds 
to 
\begin{eqnarray}
S_{l+1}(x)=x^{l+1}e^x\frac{d^l}{dx^l}\left(\frac{e^{-x}}{x}\right)
=l!L_l^{-l-1}(x),
\end{eqnarray}
while Eq.~(\ref{aux}) becomes 
\begin{eqnarray}
l L_l(x)=l L_{l-1}(x)-xL^1_{l-1}(x).
\end{eqnarray}
For Jacobi polynomials $\sigma'(x)=-2x.$ As $-2x=1-x-(1+x),$ where $1\pm x$ 
can be incorporated into the weight functions $w(x)=(1-x)^a(1+x)^b,$ there is 
no need to introduce auxiliary polynomials. For example, Eq.~(\ref{aux}) 
becomes 
\begin{eqnarray}
2l P_l^{(a,b)}(x)=(a+l)(1+x)P_{l-1}^{(a,b+1)}(x)-(b+l)(1-x)P_{l-1}^{(a+1,b)}(x).
\end{eqnarray}

\section{Discussion} 

We have used a simple and natural method for constructing polynomials 
$Q_\nu^{(\alpha,-a)}(x)$ that are complementary to the $C_n^{(\alpha,-a)}(x)$ 
polynomials and related to them by a Rodrigues formula. Similar to the 
classical orthogonal polynomials, the $Q_\nu^{(\alpha,-a)}(x)$ appear as 
solutions of a Sturm-Liouville ordinary second-order differential equation and 
obey Rodrigues formulas themselves. On the other hand, and different from the 
classical polynomials, their infinite sets of orthogonality integrals are not 
the standard ones. These real orthogonal polynomials and their nontrivial 
orthogonality properties are closely related to Romanovski polynomials and to 
physical phenomena. In summary, all basic properties of Romanovski polynomials 
derive from the Rodrigues formula~(\ref{rds}) except for the orthogonality 
integrals. 

\section{Acknowledgments}

It is a pleasure to thank M. Kirchbach for introducing me to the 
$C_n^{(\alpha,-a)}(x)$ polynomials. Thanks are also due to V. Celli for help 
with some of the orthogonality integrals.


\end{document}